\documentclass[12pt]{amsart}
\usepackage{latexsym,amssymb,amsmath}
 2

\newtheorem{thm}{Theorem}[section]
\newtheorem{lem}[thm]{Lemma}

\newtheorem{cor}[thm]{Corollary}
\newcommand{\thmref}[1]{Theorem~\ref{#1}}

\theoremstyle{remark}
\newtheorem{rmk}{Remark}[section]

\begin{document}

\title[Identities for the Ramanujan's Tau Function]
{Rankin-Cohen Brackets and  van der Pol-Type Identities for the Ramanujan's 
Tau Function} 
\author{B. Ramakrishnan and Brundaban Sahu}
\address{Harish-Chandra Research Institute, 
       Chhatnag Road, Jhunsi,
     Allahabad \linebreak
    211 019,
   India.}
\email[B. Ramakrishnan]{ramki@hri.res.in}
\email[Brundaban Sahu]{sahu@hri.res.in}


\date{\today}

\begin{abstract}
We use Rankin-Cohen brackets for modular forms and quasimodular forms to 
give a different proof of the results obtained by D. Lanphier \cite{lanphier} 
and D. Niebur \cite{niebur} on the van der Pol type identities for the 
Ramanujan's tau function. As consequences we obtain convolution sums and 
congruence relations involving the divisor functions.
\end{abstract}

\maketitle

\bigskip

\bigskip

\section{Introduction}
Let 
$$
\Delta (z) = q\prod_{n\ge 1} (1-q^n)^{24} = \sum_{n\ge 1} \tau(n) q^n
$$
be the Ramanujan's cusp form of weight $12$ for $SL_{2}({\mathbb Z})$,
where $q= e^{2\pi iz}$, $z$ belongs to the
upper half-plane ${\mathcal H}$. The function $\tau(n)$ 
is called the Ramanujan's tau-function. Using differential equations satisfied by 
$\Delta(z)$, B. van der Pol \cite{vander} derived identities relating $\tau(n)$ 
to sum-of-divisors functions. For example,

\begin{equation}\label{vander}
\tau(n) = n^2\sigma_3(n) + 60 \sum_{m=1}^{n-1} 
(2n- 3m)(n- 3m)\sigma_3(m)\sigma_3(n-m),
\end{equation}
where $\sigma_k(n) = \displaystyle{\sum_{d|n}}d^k$. Using the relation between 
$\sigma_3(n)$ and $\sigma_7(n)$ (see Theorem 3.1 (i) below), this is equivalent 
to the following identity: 

\begin{equation}\label{vander1}
\tau(n) = n^2\sigma_7(n) - 540 \sum_{m=1}^{n-1} 
m(n-m)\sigma_3(m)\sigma_3(n-m).
\end{equation}

\noindent In \cite{lanphier}, D. Lanphier used differential operators studied by 
Maass \cite{maass} to prove the above van der Pol identity \eqref{vander1}. 
He also obtained several van der Pol-type identities using the Maass operators 
and thereby obtained new congruences for the Ramanujan's tau-function. 

In \cite{niebur}, D. Niebur derived a formula for $\tau(n)$ similar to 
the classical ones of Ramanujan and van der Pol, but has the feature that 
higher divisor sums do not appear (see \thmref{niebur:thm} (i) below). 
Niebur proved the formula by expressing 
$\Delta(z)$ in terms of the logarithmic derivatives of $\Delta(z)$. 

In this paper we show that the identities for the Ramanujan function 
$\tau(n)$ proved by Lanphier in \cite{lanphier} can be obtained using 
the Rankin-Cohen bracket for modular forms and the basic relations among 
the Eisenstein series.  
By our method we also obtain new identities for $\tau(n)$ which are not proved 
in \cite{lanphier}. Next we show that the theory of quasimodular forms 
can be used to prove Niebur's formula for $\tau(n)$. The method of using 
quasimodular forms gives new formulas for $\tau(n)$. Though one can obtain many 
identities for $\tau(n)$, here we restrict only those identities in which 
only the convolution of divisor functions appear. Finally, we observe that  
the identities of $\tau(n)$ give rise to various identities for the 
convolution of the divisor functions. As a consequence, we also present some 
congruences involving the divisor functions.

\section{Statement of results}
\begin{thm}
\begin{eqnarray*}
{\rm (i)}\quad \tau(n) &=& n^2\sigma_7(n) - 540 
\sum_{m=1}^{n-1}m(n-m)\sigma_3(m)\sigma_3(n-m), \\
{\rm (ii)}\quad \tau(n) &=& -\frac{5}{4} n^2\sigma_7(n) + \frac{9}{4} 
n^2\sigma_3(n)
+ 540\sum_{m=1}^{n-1}m^2\sigma_3(m)\sigma_3(n-m), \\
{\rm (iii)}\quad \tau(n) &=& n^2\sigma_7(n) - \frac{1080}{n}
\sum_{m=1}^{n-1}m^2(n-m)\sigma_3(m)\sigma_3(n-m), \\
{\rm (iv)}\quad \tau(n) &=& -\frac{1}{2} n^2\sigma_7(n) + \frac{3}{2} 
n^2\sigma_3(n)
+ \frac{360}{n}\sum_{m=1}^{n-1}m^3\sigma_3(m)\sigma_3(n-m). 
\end{eqnarray*}
\end{thm}
\begin{thm}   
\begin{eqnarray*}
{\rm (i)}\quad \tau(n) &=& -\frac{11}{24} n\sigma_9(n) + \frac{35}{24} 
n\sigma_5(n) + {350}\sum_{m=1}^{n-1}(n-m)\sigma_3(m)\sigma_5(n-m), \\
{\rm (ii)}\quad \tau(n) &=& \frac{11}{36} n\sigma_9(n) + \frac{25}{36} 
n\sigma_3(n)- {350}\sum_{m=1}^{n-1}m\sigma_3(m)\sigma_5(n-m), \\
{\rm (iii)}\quad \tau(n) &=& \frac{1}{6} n\sigma_9(n) + \frac{5}{6} n\sigma_3(n)
- \frac{420}{n}\sum_{m=1}^{n-1}m^2\sigma_3(m)\sigma_5(n-m), \\
{\rm (iv)}\quad \tau(n) &=&  n\sigma_9(n)  
- \frac{2100}{n}\sum_{m=1}^{n-1}m(n-m)\sigma_3(m)\sigma_5(n-m), \\
{\rm (v)}\quad \tau(n) &=& -\frac{1}{4} n\sigma_9(n) + \frac{5}{4} n\sigma_5(n)
+ \frac{300}{n}\sum_{m=1}^{n-1}(n-m)^2\sigma_3(m)\sigma_5(n-m). 
\end{eqnarray*}
\end{thm}

\begin{thm}
\begin{eqnarray*}
\tau(n) &=& \frac{65}{756} \sigma_{11}(n) + \frac{691}{756} \sigma_5(n)
- \frac{2\cdot 691}{3n}\sum_{m=1}^{n-1}m\sigma_5(m)\sigma_5(n-m). 
\end{eqnarray*}
\end{thm}

\begin{thm}
\begin{eqnarray*}
{\rm (i)}\quad \tau(n) &=& -\frac{91}{600} \sigma_{11}(n) + \frac{691}{600} 
\sigma_3(n)
+ \frac{4\cdot 691}{5n}\sum_{m=1}^{n-1}m\sigma_3(m)\sigma_7(n-m),\\   
{\rm (ii)}\quad \tau(n) &=& -\frac{91}{600} \sigma_{11}(n) + \frac{691}{600} 
\sigma_7(n)
+ \frac{2\cdot 691}{5n}\sum_{m=1}^{n-1}(n-m)\sigma_3(m)\sigma_7(n-m).
\end{eqnarray*}
\end{thm}

\begin{rmk}
Theorems 2.1 to 2.4 are exactly the same as Theorems 1 to 4 of 
\cite{lanphier}.
\end{rmk}

The following theorem gives identities for $\tau(n)$ in which the convolution 
part contains only the divisor function $\sigma(n) (=\sigma_1(n))$. 

\begin{thm}\label{niebur:thm}
\begin{eqnarray*}
{\rm (i)}\quad \tau(n) &=&
n^4\sigma(n) - 24\sum_{m=1}^{n-1} (35 m^4 -52 m^3n + 18 m^2 n^2) 
\sigma(m)\sigma(n-m), \\
{\rm (ii)}\quad \tau(n) &=& n^4(7\sigma(n) - 6 \sigma_3(n)) 
- 168\sum_{m=1}^{n-1} (5 m^4 - 4 m^3n)
\sigma(m)\sigma(n-m), \\
{\rm (iii)}\quad \tau(n) &=& n^4\sigma_3(n) - 
 168\sum_{m=1}^{n-1} (5 m^4 - 8 m^3n + 3 m^2n^2)
\sigma(m)\sigma(n-m), \\
{\rm (iv)}\quad \tau(n) &=& \frac{n^4}{3}(7\sigma(n) -4\sigma_3(n)) - 
56\sum_{m=1}^{n-1} (15 m^4 - 20 m^3n + 6m^2n^2)
\sigma(m)\sigma(n-m). 
\end{eqnarray*}
\end{thm}

\begin{rmk}
In Theorem 2.5, formula (i) was proved by Niebur \cite{niebur}. 
The rest of the formulas are obtained while proving (i) by using 
quasimodular forms. 
\end{rmk}

In the above theorems we presented results that are already proved by 
Niebur and Lanphier and a few more identities like Niebur's.
We now state some more new identities  using our method. 
As mentioned in the introduction we restrict ourselves to the identities 
which involve only the convolution of the divisor functions. We list these 
identities in two theorems, one uses the theory of modular forms and the 
other uses the theory of quasimodular forms.

\begin{thm}
\begin{eqnarray*}
{\rm (i)}\quad \tau(n) &=& \frac{5}{12} n\sigma_{3}(n) + \frac{7}{12} 
n\sigma_5(n)
+70 \sum_{m=1}^{n-1}(2n-5m)\sigma_3(m)\sigma_5(n-m), \\
{\rm (ii)}\quad \tau(n) &=& n^2\sigma_{3}(n) + 
60 \sum_{m=1}^{n-1}(4n^2-13mn+9m^2)\sigma_3(m)\sigma_3(n-m). \\
\end{eqnarray*}
\begin{eqnarray*}
{\rm (iii)}\quad \tau(n) &=& \frac{65}{756} \sigma_{11}(n) + 
\frac{5\cdot691}{12\cdot 756 n} \sigma_7(n) 
+ \frac{691}{12\cdot 108 n} \sigma_5(n)  \\
& & \hskip 3cm   - \frac{5\cdot 691}{54 n^2} \sum_{m=1}^{n-1}(3n - 7m)
\sigma_5(m)\sigma_7(n-m) \\
& & \hskip 3cm  -\frac{13\cdot 691}{9 n^2}\sum_{m=1}^{n-1}m(n-m)
\sigma_5(m)\sigma_5(n-m),\\
{\rm (iv)}\quad \tau(n) &=& \frac{65}{756} \sigma_{11}(n) + 
\frac{3\cdot691}{8\cdot441} \sigma_3(n) + \frac{5\cdot 691}{24\cdot441} 
\sigma_7(n)\\
& & \qquad + \frac{5\cdot 691}{441 n^2} 
\sum_{m=1}^{n-1}(91m^2 - 65 mn + 10n^2) \sigma_3(m)\sigma_7(n-m) \nonumber 
\\
& & \qquad \qquad -\frac{13\cdot 691}{9 n^2}\sum_{m=1}^{n-1}m(n-m)
\sigma_5(m)\sigma_5(n-m),\\
{\rm (v)}\quad \tau(n) &=& \frac{65}{756} \sigma_{11}(n) + 
\frac{25\cdot691}{36\cdot 756 n} \sigma_3(n) 
+ \frac{11\cdot691}{36\cdot 756 n} \sigma_9(n) \nonumber \\
& & \hskip 3cm   - \frac{55\cdot 691}{1134 n^2} \sum_{m=1}^{n-1}(7m-2n)
\sigma_3(m)\sigma_9(n-m) \\
& & \hskip 3cm  -\frac{13\cdot 691}{9 n^2}\sum_{m=1}^{n-1}m(n-m)
\sigma_5(m)\sigma_5(n-m).\\
\end{eqnarray*}
\end{thm}

\begin{thm}
\begin{eqnarray*}
{\rm (i)} \quad \tau(n) &=& \frac{5\cdot 691}{9504}\sigma(n) - 
\frac{(6n-5)\cdot 691}{864} \sigma_9(n) \\
& & \qquad \quad + \frac{2275}{1584}\sigma_{11}(n)
- \frac{5\cdot 691}{864} \sum^{n-1}_{m=1}\sigma(m)\sigma_9(n-m),\\
{\rm (ii)} \quad \tau(n) &=&
\frac{15}{32}n\sigma(n)-\frac{33}{32}n\sigma_9(n)+\frac{50}{32}n^2\sigma_7(n)
+225 \sum^{n-1}_{m=1}m\sigma(m)\sigma_7(n-m),\\
{\rm (iii)} \quad \tau(n) &=&
\frac{6}{7}n^2\sigma(n) -\frac{9}{7}n^3\sigma_5(n) +\frac{10}{7}n^2\sigma_7(n)
-24\cdot 18 \sum^{n-1}_{m=1}m^2\sigma(m)\sigma_5(n-m),\\
{\rm (iv)} \quad \tau(n) &=&
 \frac{14}{5}n^3\sigma(n)+\frac{12}{5}n^4\sigma_3(n)-\frac{21}{5}n^3\sigma_5(n)
+ 24\cdot 28 \sum^{n-1}_{m=1}m^3\sigma(m)\sigma_3(n-m), \\
{\rm (v)} \quad \tau(n) &=&
\frac{5}{12}n\sigma(n)+\frac{25}{24}n\sigma_7(n)-\frac{11}{24}n\sigma_9(n)
+25\sum^{n-1}_{m=1}(9m-n)\sigma(m)\sigma_7(n-m),\\
\end{eqnarray*}
\begin{eqnarray*}
{\rm (vi)} \quad \tau(n) &=&
\frac{9}{14}n^2\sigma(n)+\frac{5}{14}n^2\sigma_7(n)
-108 \sum^{n-1}_{m=1}(4m^2-mn)\sigma(m)\sigma_5(n-m),\\
{\rm (vii)} \quad \tau(n) &=& 
\frac{8}{5}n^3\sigma(n)-\frac{3}{5}n^3\sigma_5(n)
+96 \sum^{n-1}_{m=1}(7m^3-3m^2n)\sigma(m)\sigma_3(n-m),\\
{\rm (viii)} \quad \tau(n) &=& \frac{1}{2}n^2\sigma(n) + 
\frac{1}{2}n^2\sigma_5(n) - 
12 \sum^{n-1}_{m=1}(36m^2-16 mn +n^2)\sigma(m)\sigma_5(n-m),\\
{\rm (ix)} \quad \tau(n) &=&
n^3\sigma(n)- 24 \sum^{n-1}_{m=1}(21m^2n-28m^3-3mn^2)\sigma(m)\sigma_3
(n-m).
\end{eqnarray*}
\end{thm}

As a consequence of the above theorems, we get the following congruences 
for $\tau(n)$ in terms of the divisor functions. Some of the congruences 
are already known (for example (iii) and (iv) are corresponding to 
(7.15) and (5.6) of \cite{lahiri2}) and some are new. 

\begin{cor}\label{tau-cong}
\begin{eqnarray*}
{\rm (i)} \qquad 12\tau(n) &\equiv& 5n\sigma_3(n)+7n\sigma_5(n) \pmod
{2^3\cdot 3\cdot 5\cdot 7},\\
{\rm (ii)} \qquad 32\tau(n) &\equiv&
15n\sigma(n)+ 50n^2\sigma_7(n)-33n\sigma_9(n) \pmod {2^5\cdot 3^2\cdot 5^2},\\
{\rm (iii)} \qquad 7\tau(n) &\equiv& 6n^2\sigma(n)- 9n^3\sigma_5(n)
+ 10 n^2\sigma_7(n) \pmod{2^4\cdot 3^3\cdot 7},\\
{\rm (iv)} \qquad 5\tau(n) &\equiv&
14 n^3\sigma(n)+ 12 n^4\sigma_3(n)- 21 n^3\sigma_5(n) \pmod {2^5\cdot 3\cdot 
5\cdot 7},\\
{\rm (v)} \qquad 24\tau(n) &\equiv& 10 n\sigma(n)+ 25 n\sigma_7(n)- 11 
n\sigma_9(n) \pmod {2^3\cdot 3\cdot 5^2},\\
{\rm (vi)} \qquad 14\tau(n) &\equiv& 9 n^2\sigma(n)+ 5 n^2\sigma_7(n)  
\pmod {2^3\cdot 3^3\cdot 7},\\
{\rm (vii)} \qquad 5\tau(n) &\equiv& 8 n^3\sigma(n)- 3 n^3\sigma_5(n)
\pmod {2^5\cdot 3\cdot 5},\\
{\rm (viii)} \qquad 2\tau(n) &\equiv& n^2\sigma(n)+ n^2\sigma_5(n)
\pmod {2^3\cdot 3}.
\end{eqnarray*}
\end{cor}

\bigskip

Finally, we state some identities for certain convolution
of divisor functions. Some of these identities are obtained while 
using our method and the rest follow from the identities for $\tau(n)$ 
mentioned in the above theorems. We remark that except for (vii), all other 
formulas are different from the ones obtained by E. Royer \cite{royer} .

\begin{thm}\label{sigma:identities}
\begin{eqnarray*}
{\rm (i)} \qquad \sum_{m=1}^{n-1}m^3 \sigma(m)\sigma(n-m) &=& 
\frac{1}{12}n^3\sigma_3(n) - \frac{1}{24}n^3(3n-1)\sigma(n),\\
{\rm (ii)}\qquad   \sum_{m=1}^{n-1}m^2 \sigma(m)\sigma(n-m) &=& 
\frac{1}{8}n^2\sigma_3(n) - \frac{1}{24}n^2(4n-1)\sigma(n), \\
\end{eqnarray*}
\begin{eqnarray*}
{\rm (iii)}~~\quad   \sum_{m=1}^{n-1}m \sigma(m)\sigma(n-m) &=& 
\frac{1}{24}n(1 -6n)\sigma(n) + \frac{5}{24}n\sigma_3(n),\\
{\rm (iv)} \quad \sum^{n-1}_{m=1}m^2\sigma(m)\sigma_3(n-m) & = & 
- \frac{1}{240}n^2\sigma(n)-\frac{1}{120}n^2\sigma_3(n)+ 
\frac{1}{80}n^2\sigma_5(n),  \\
{\rm (v)}~~ \quad \sum^{n-1}_{m=1}m\sigma(m)\sigma_3(n-m) & = & 
- \frac{1}{240}n\sigma(n)-\frac{1}{40}n^2\sigma_3(n)+ 
\frac{7}{240}n\sigma_5(n),  \\
{\rm (vi)}~~ \quad \sum^{n-1}_{m=1}m\sigma(m)\sigma_5(n-m) & = & 
\frac{1}{504}n\sigma(n) - \frac{1}{84}n^2\sigma_5(n) + 
\frac{5}{504}n\sigma_7(n), \\
{\rm (vii)} \qquad \sum^{n-1}_{m=1}\sigma(m)\sigma_5(n-m) & = & 
\frac{1}{504}\sigma(n)-\frac{1}{12}n\sigma_5(n) + \frac{1}{24} \sigma_5(n) + 
\frac{5}{126}\sigma_7(n), \\
 {\rm (viii)} \qquad \sum^{n-1}_{m=1}\sigma(m)\sigma_7(n-m)&=& 
-\frac{1}{480}\sigma(n)+\frac{1}{24} \sigma_7(n)+\frac{11}{480}\sigma_9(n)
- \frac{1}{16}n\sigma_7(n). \\
\end{eqnarray*}
\end{thm}

\begin{cor}
\begin{equation}
\sum_{m=1}^{n-1} (2m^3 - 3 m^2 n + m n^2) \sigma(m) \sigma(n-m) = 0.
\end{equation}
\end{cor}

Using Theorem 2.8 (i)--(iii) and Theorem 2.5 (i), we get another formula for 
$\tau(n)$, given in the following corollary. 

\begin{cor}
\begin{equation}
\tau(n) = 50 n^4 \sigma_3(n) - 7 n^4(12 n - 5)\sigma(n) - 840\sum_{m=1}^{n-1} 
m^4 \sigma(m)\sigma(n-m). 
\end{equation}
\end{cor}

We end by stating some  congruence relations among the divisor functions. 
These congruences follow as a consequence of the above convolution identities 
and the congruences of $\tau(n)$ (Corollary 2.8). 

\begin{cor}
\begin{eqnarray*}
\!\!{\rm (i)} ~~\quad\qquad\qquad (6n-5) \sigma(n) & \equiv & \sigma_3(n) 
\pmod{24}, ~~ \gcd(n,6)=1, \\
\!\!{\rm (ii)}~~~~ \qquad \qquad \sigma(n)+2n\sigma_3(n) & \equiv & 3 
\sigma_5(n) \pmod{16}, ~~2\not\vert n,\\
{\rm (iii)} \qquad \qquad  n\sigma(n)+5n\sigma_7(n) & \equiv& 
6n^2\sigma_5(n) \pmod{2^3\cdot 3^2\cdot 7}, ~~\gcd(n, 42)=1,\\
{\rm (iv)} ~~\qquad \quad  20 \sigma_7(n)+11\sigma_9(n) &\equiv& 
\sigma(n)+30n\sigma_7(n) \pmod{2^5\cdot 3\cdot5},\\
{\rm (v)} \qquad \qquad  5 \sigma(n)+6n\sigma_7(n) &\equiv& 
11\sigma_9(n) \pmod{2^5}, ~~2\not\vert n,\\
{\rm (vi)} ~~~\qquad \qquad  \sigma(n)+2n\sigma_3(n) &\equiv& 
3\sigma_5(n) \pmod{2^4\cdot5}, ~~\gcd(n,10)=1,\\
{\rm (vii)} \quad   \sigma(n)+10(3n-2)\sigma_7(n) &\equiv& 
11\sigma_9(n) \pmod{2^3\cdot3\cdot5}, ~~\gcd(n, 30)=1.
\end{eqnarray*}
\end{cor}

\section{Preliminaries}

In this section we shall provide some well-known facts about 
modular forms and give the definitions of Rankin-Cohen bracket and 
quasimodular forms, which are essential in proving our theorems. For 
basic details of the theory of modular forms and quasimodular forms, we 
refer to \cite{{serre},{kaneko},{martin-royer1},{royer}}. 

For an even integer $k\ge 4$, let $M_k(1)$ (resp. $S_k(1)$) denote the 
vector space of modular forms (resp. cusp forms) of weight $k$ for the 
full modular group $SL_2({\mathbb Z})$. Let $E_k$ be the normalized 
Eisenstein series of weight $k$ in $M_k(1)$, given by 
$$
E_k(z) = 1 - \frac{4k}{B_k}\sum_{n\ge 1} \sigma_{k-1}(n) q^n,
$$
where $B_k$ is the $k$-th Bernoulli number defined by 
$$
\frac{x}{e^x-1} = \sum_{m=0}^\infty \frac{B_m}{m!} x^m.
$$
The first few Eisenstein series are as follows.

\begin{equation}
\begin{split}
E_4 & = 1 + 240\sum_{n\ge 1}\sigma_3(n) q^n, \\
E_6 & = 1 - 504\sum_{n\ge 1}\sigma_5(n) q^n, \\
E_8 & = 1 + 480\sum_{n\ge 1}\sigma_7(n) q^n, \\
E_{10} & = 1 - 264\sum_{n\ge 1}\sigma_9(n) q^n, \\
E_{12} & = 1 + \frac{65520}{691}\sum_{n\ge 1}\sigma_{11}(n) q^n. 
\end{split}
\end{equation}

Since the space $M_k(1)$ is one-dimensional for $4\le k \le 10$ and for  
$k=14$ and $S_{12}(1)$ is one-dimensional, we have the following 
well-known identities:
\begin{thm}$~~$\\

\begin{center}
\begin{tabular}{rrcl}
{\rm (i)}\qquad  & $E_8(z)$  & = & $E_4^2(z)$,\hskip 2cm  \\
&&&\\
{\rm (ii)} \quad &  $E_{10}(z)$ & = & $E_4(z) E_6(z)$,\hskip 2cm \\
&&&\\
{\rm (iii)}\qquad & $E_{12}(z) - E_8(z) E_4(z)$  & = & $\left(\frac{65520}{691} 
- 720\right) \Delta(z)$, \hskip 2cm \\
&&&\\
{\rm (iv)}\qquad  & $E_{12}(z) - E_6^2(z)$  & = & $\left(\frac{65520}{691} + 
1008\right) \Delta(z)$,\hskip 2cm 
\end{tabular}
\end{center}

\end{thm}
  
\bigskip

\noindent {\bf Rankin-Cohen Brackets}: The derivative of a modular form is not a 
modular form. However, the works of Rankin and Cohen \cite{{rankin},{cohen}} 
lead to the concept of Rankin-Cohen brackets which is defined in the following. 

\smallskip

\noindent {\bf Definition}:~ 
Let $f \in M_k(1)$ and $g\in M_l(1)$ be modular forms of weights $k$ and $l$ 
respectively.  For each $\nu \ge 0,$ define the $\nu$-th {\em Rankin-Cohen 
bracket} (in short RC bracket) of $f$ and $g$ by
\begin{equation}
[f, g]_{\nu} := \sum^\nu_{r=0}(-1)^r \binom{\nu + k-1}{\nu-r} 
\binom{\nu+l-1}{r} D^{(r)} f~ D^{(\nu -r)}g,
\end{equation}
where we have set $D:= \displaystyle{\frac{1}{2\pi i} \frac{d}{dz}}$. When 
$\nu=1$, we write $[f,g]$ instead of $[f,g]_1$. 

\begin{thm} (\cite[pp. 58--61]{zagier})
Let $f \in M_k(1)$ and $g\in M_l(1)$. Then $[f,g]_\nu$ is a modular form of 
weight $k+l+2\nu$ for $SL_2({\mathbb Z})$. It is a cusp form if $\nu \ge 1$.
\end{thm}

\bigskip

\noindent {\bf Quasimodular forms}: We now present some basics of quasimodular 
forms. Another important Eisenstein series is the weight $2$ Eisenstein series 
$E_2$ given by 
\begin{equation}
E_2(z) = 1 - 24 \sum_{n\ge 1}\sigma(n) q^n.
\end{equation}
It is not a modular form because it doesn't satisfy the required 
transformation property under the action of $SL_2({\mathbb Z})$. 
However, it plays a fundamental role in defining the concept of quasimodular 
forms, which was formally introduced by M. Kaneko and D. Zagier \cite{kaneko}. 

\smallskip

\noindent {\bf Definition}:~ Let $k\ge 1, s\ge 0$ be natural numbers.
A holomorphic function $f:{\mathcal H} 
\rightarrow {\mathbb C}$ is defined to be a quasimodular form of weight 
$k$, depth $s$ on $SL_2({\mathbb Z})$, if there exist holomorphic functions 
$f_0, f_1, \ldots, f_s$ on ${\mathcal H}$ such that 
\begin{equation}
(cz+d)^{-k} f\left(\frac{az+b}{cz+d}\right) = \sum_{i=0}^{s}f_i(z) 
\left(\frac{c}{cz+d}\right)^i,
\end{equation}
for all $\begin{pmatrix}a&b \\ c&d\\\end{pmatrix} \in SL_2({\bf Z})$ and 
such that $f_s$ is holomorphic at infinity and not identically vanishing. 

\smallskip

\begin{rmk}
It is a fact that if $f$ is a quasimodular form of weight $k$ and depth 
$s$, not identically zero, then $k$ is even and $s\le k/2$. 
\end{rmk}

\smallskip

The space of all quasimodular forms of weight $k$, depth $s$ on $SL_2({\mathbb 
Z})$ is denoted by $\widetilde{M}_k^{\le s}(1)$. Note that $E_2$ is a 
quasimodular form of weight $2$ and depth $1$, and so $E_2\in 
\widetilde{M}_2^{\le 1}(1)$. 

\smallskip

We need the following lemma (see \cite{royer} for details).

\begin{lem}(\cite[Lemma 1.17]{royer})~ Let $k\ge 2$ be even. Then  \\
$$
\widetilde{M}_k^{\le k/2}(1) = \bigoplus_{i=0}^{k/2-1} 
D^iM_{k-2i}(1) \oplus {\mathbb C} D^{k/2-1}E_2.
$$
\end{lem}

\smallskip

\noindent {\bf Definition}:~(Rankin-Cohen bracket of quasimodular forms)~ 
Let $f \in \widetilde{M}_k^{\le s}(1)$ and $g\in \widetilde{M}_l^{\le t}(1)$ 
be quasimodular forms of weights $k$, $l$ and depths $s$, $t$  
respectively.  For each $n \ge 0,$ define the $n$-th {\em Rankin-Cohen 
Bracket} (in short RC bracket) of $f$ and $g$ by
\begin{equation}
\Phi_{n;k,s;l,t}(f,g):= \sum^\nu_{r=0}(-1)^r \binom{k-s + n-1}{n-r} 
\binom{l-t+n-1}{r} D^{(r)} f~ D^{(n-r)}g.
\end{equation}
\begin{thm}(\cite[Theorem 1]{martin-royer2})~Let $k,l \ge 2$ and $s,t\ge 0$ with 
$s\le k/2$ and $t\le l/2$. Then for $n\ge 0$, the RC bracket 
$\Phi_{n;k,s;l,t}(f,g)$ is a quasimodular form of weight $k+l+2n$ and depth 
$s+t$ for $SL_2({\mathbb Z})$. 
\end{thm}

\section{proofs}

\subsection{Proof of Theorem 2.1:}
Differentiating (i) of Theorem 3.1 twice, we get 
\begin{equation}\label{1}
D^2E_8 = 2(DE_4)^2 + 2E_4 D^2E_4. 
\end{equation}
Now consider the RC bracket $[E_4,E_4]_2 = 20 D^2E_4 E_4 - 
25 (DE_4)^2 \in S_{12}(1) = {\mathbb C}\Delta $. Substituting for 
$E_4 D^2E_4$ from \eqref{1}, we get 
\begin{equation}\label{2}
2 D^2E_8 - 9 (DE_4)^2  = 960 \Delta.
\end{equation}
Comparing the $n$-th Fourier coefficients we obtain (i). 

On the other hand, in \eqref{2} substituting for $(DE_4)^2$ from 
\eqref{1}, we get 
$$
-5 D^2E_8 + 18 E_4 D^2E_4 = 1920 \Delta,
$$
from which we obtain the identity (ii). To obtain (iii) we differentiate 
\eqref{2} and compare the coefficients. Finally to prove (iv), we differentiate 
\eqref{1} and use it in \eqref{2} to get 
$$
960 D\Delta = -D^3E_8 + 6E_4D^3E_4. 
$$
Now comparing the $n$-th Fourier coefficients we get (iv). This completes the 
proof of Theorem 2.1. 

\subsection{Proof of Theorem 2.2:} We now take differentiation of Theorem 
3.1 (ii) to get 
\begin{equation}\label{3}
DE_{10} = E_6DE_4  + E_4 DE_6.
\end{equation}
Now consider the RC bracket $[E_4, E_6]$ which belongs to $S_{12}(1)$ and so 
is a constant multiple of $\Delta$. Thus, we have 
\begin{equation}\label{4}
4E_4 DE_6 - 6 E_6 DE_4 = -3456~\Delta.
\end{equation}
Substituting for $E_4 DE_6$ from \eqref{3} in \eqref{4} and comparing the 
$n$-th Fourier coefficients, we get the identity (i). Instead, if we 
substitute for $E_6 DE_4$ from \eqref{3} in \eqref{4}, we get the identity 
(ii). Taking derivative of \eqref{3} we get 
\begin{equation}\label{5}
D^2E_{10} = E_6D^2E_4  + 2 (DE_4)(DE_6) + E_4 D^2E_6.
\end{equation}
Differentiating \eqref{4} we get 
\begin{equation}\label{6}
-2 (DE_4)(DE_6) + 4 E_4 D^2E_6 - 6 E_6D^2E_4 = -3456~D\Delta.
\end{equation}
Eliminating $E_4 D^2E_6$ from \eqref{5} and \eqref{6} we obtain
\begin{equation}\label{7}
2 D^2E_{10} - 5 (E_6D^2E_4 + (DE_4)(DE_6)) = -1728~D\Delta.
\end{equation}
Now, consider the RC bracket $[E_4,E_6]_2$. Since it belongs 
to $S_{14}(1)$, it must be zero. So, we get the following 
\begin{equation}\label{8}
10 E_4 D^2E_6 - 35 (DE_4)(DE_6) + 21 E_6D^2E_4 = 0.
\end{equation}
Eliminating $(DE_4)(DE_6)$ from \eqref{5} and \eqref{7}, we get 
\begin{equation}\label{9}
- D^2E_{10} - 5 E_6 D^2E_4 + 5 E_4 D^2E_6 + 3456~D\Delta = 0, 
\end{equation}
and from \eqref{5} and \eqref{8}, we get 
\begin{equation}\label{10}
-35 D^2E_{10} + 77 E_6 D^2E_4 + 55 E_4 D^2E_6 = 0.
\end{equation}
Now, eliminating $E_4 D^2E_6$ from \eqref{9} and \eqref{10}, we obtain 
\begin{equation*}
24 D^2E_{10} - 132 E_6 D^2E_4 + 38016~D\Delta = 0,
\end{equation*}
whose $n$-th Fourier coefficient give the identity (iii). In the last step 
if we eliminate $E_6 D^2E_4$ instead of $E_4 D^2E_6$, we obtain the identity 
(v). For the proof of identity (iv), we first eliminate $E_6 D^2E_4$ from 
\eqref{5} \& \eqref{7} and \eqref{5} \& \eqref{8} to get 
\begin{eqnarray*}
-3 D^2E_{10} + 5 (DE_4)(DE_6) + 5 E_4 D^2E_6 + 1728~D\Delta &= &0, \\
-21 D^2E_{10} + 77 (DE_4)(DE_6) + 11 E_4 D^2E_6 &= &0. 
\end{eqnarray*}
The required identity (iv) follows by eliminating $E_4 D^2E_6$ from the above 
two equations. This completes the proof of Theorem 2.2.  

\smallskip

\subsection{Proof of Theorem 2.3:} This identity follows easily by 
differentiating (iv) of Theorem 3.1 and comparing the $n$-th Fourier 
coefficients. 

\smallskip

\subsection{Proof of Theorem 2.4:} Differentiating (iii) of Theorem 3.1 
gives 
\begin{equation}\label{11}
DE_{12} - E_4 DE_8 - E_8 DE_4 = \left(\frac{65520}{691} -
720\right) D\Delta.
\end{equation}
Next, using the fact that the RC bracket $[E_4, E_8]$ is zero (because 
it belongs to the space $S_{14}(1)$) we get 
\begin{equation}\label{12}
E_4DE_8 = 2 E_8 DE_4.
\end{equation}
Substituting  for $E_4DE_8$ from \eqref{12} in \eqref{11} and 
comparing the $n$-th Fourier coefficients yield identity (i). On the other 
hand, substituting for $E_8DE_4$ from \eqref{12} in \eqref{11} 
gives 
\begin{equation}
DE_{12} - \frac{3}{2} E_4 DE_8 = \left(\frac{65520}{691} -
720\right) D\Delta,
\end{equation}
from which (ii) follows. 

\smallskip

\subsection{Proof of Theorem 2.5:} As mentioned before, we make use 
of the theory of quasimodular forms. Recall that the Eisenstein series $E_2$ is a 
quasimodular form of weight $2$ and  depth $1$ on $SL_2({\mathbb Z})$. So, 
$D^iE_2 \in \widetilde{M}_{2i+2}^{\le i+1}(1)$. Consider the following $6$ 
quasimodular forms of weight $12$ on $SL_2({\mathbb Z})$ which are 
RC brackets of functions involving $D^iE_2$.\\

\begin{center}
\begin{tabular}{rll}
{\rm (i)}\quad $f_1(z)$ &=~ $\Phi_{1;8,4;2,1}(D^3E_2, E_2)$ &=~
$4 D^3 E_2 DE_2- E_2 D^4 E_2$,\\
&&\\
{\rm (ii)}\quad $f_2(z)$ &=~  $ \Phi_{1;6,3;4,2}(D^2E_2, D E_2)$ &=~ 
$3 (D^2 E_2)^2 - 2 D^3 E_2 D E_2$,\\
&&\\
{\rm (iii)}\quad $f_3(z)$ &=~  $\Phi_{2;6,3;2,1}(D^2E_2, E_2)$ &=~  
$6 (D^2 E_2)^2 - 8 D^3E_2 D E_2 + E_2 D^4 E_2$,\\
&&\\
{\rm (iv)}\quad $f_4(z)$ &=~  $\Phi_{2;4,2;4,2}(D E_2, D E_2)$ &=~  
$6 D^3 E_2 DE_2- 9 (D^2E_2)^2$,\\
&&\\
{\rm (v)}\quad $f_5(z)$ &=~ $\Phi_{3;4,2;2,1}(D E_2, E_2)$  &=~  
$16 D^3 E_2 DE_2- 18 (D^2E_2)^2 -  E_2 D^4 E_2$,\\
&&\\
{\rm (vi)}\quad $f_6(z)$ &=~ $\Phi_{4;2,1;2,1}(E_2, E_2)$ &=~  
$-32 D^3 E_2 DE_2 +36 (D^2E_2)^2 +2 E_2 D^4E_2$. \\
&&\\
\end{tabular}
\end{center}

\noindent 
First note that  $f_6(z) = -2 f_5(z)$ and $f_4(z) = -3 f_2(z)$.  
Using the definition of the RC 
brackets one finds that $f_1,f_2 \in  \widetilde{M}_{12}^{\le 5}(1)$, $f_3,f_4 
\in \widetilde{M}_{12}^{\le 4}(1)$, $f_5 \in \widetilde{M}_{12}^{\le 3}(1)$ 
and $f_6 \in \widetilde{M}_{12}^{\le 2}(1)$. Considering all the functions 
$f_i$, $1\le i\le 6$ in the space  $\widetilde{M}_{12}^{\le 6}(1)$ and 
using the decomposition stated in  Lemma 3.3, we have the following expressions 
for the 
$f_i$'s.

\begin{equation}\label{fi}
\begin{split}
f_1(z) &= \frac{24}{7} \Delta(z) + \frac{3}{35} D^4E_4, \\
f_2(z) &= -\frac{24}{7} \Delta(z) + \frac{1}{70} D^4E_4,\\
f_3(z) &= -\frac{72}{7} \Delta(z) - \frac{2}{35} D^4E_4,\\
f_5(z) &= 24 \Delta(z), \\
\end{split}
\end{equation}

\noindent Using the definitions of $f_i$ and \eqref{fi}, we express the 
RC brackets appearing in the definitions of $f_i$ in terms of $\Delta(z)$ 
and $D^4E_4(z)$. This way we get four expressions corresponding to the functions 
$f_1, f_2, f_3$ and $f_5$. Comparing the Fourier expansions corresponding to 
$f_5$ gives the required Niebur's identity (i). It also follows by 
eliminating $D^4E_4$ between any two expressions corresponding to $f_i$'s,
$1\le i\le 3$. The rest of the  identities ((ii) to (iv)) 
are obtained  directly by comparing the Fourier coefficients of the
expressions corresponding to $f_i$, $1\le i\le 3$, respectively.

\bigskip

\subsection{Proof of Theorem 2.6:} The identities (i) and (ii) follow directly 
by comparing the Fourier coefficients of the RC brackets 
$[E_4, E_6]_1 = 3456 \Delta$ and $[E_4, E_4]_2 = 960 \Delta$ respectively. 
Now we  prove (iii).  Differentiating  twice the expression (iv) of Theorem 3.1 
we get 
\begin{equation}\label{13}
2(DE_6)^2 + 2 E_6 D^2E_6 = D^2E_{12}- \alpha D^2\Delta,
\end{equation}
where $\alpha=\left(\frac{65520}{691}+1008\right)$. Since $\dim S_{16}(1) =1$, 
we get $[E_6, E_6]_2 = 42 D^2E_6-49(DE_6)^2 = - \frac{49}{48} [E_6, E_8] = 
- \frac{49}{48} (8 E_8 DE_6 - 6 E_6 DE_8)$. Substituting for $E_6D^2E_6$ from 
equation \eqref{13} and comparing the Fourier coefficients 
both the sides give the identity (iii).  
In the above argument, replacing $[E_6, E_8]$ by $[E_4, E_8]_2$ one gets the 
identity (iv) and replacing by $[E_4, E_{10}]$ one gets (v). 

\bigskip

\subsection{Proof of Theorem 2.7:} All the identities are straight forward 
comparison of the Fourier coefficients of the following quasimodular forms of 
weight $12$ (using the decomposition given in Lemma 3.3) respectively:
$E_2E_{10}, DE_2E_8$,  $D^2E_2E_6$, $D^3E_2E_4$, $\Phi_{1;2,1;8,0}(E_2, E_8)$, 
$\Phi_{1;4,2;6,0}(DE_2, E_6)$, $\Phi_{1;6,3;4,0}(D^2E_2, E_4)$,
$\Phi_{2;2,1;6,0}(E_2, E_6)$ and $\Phi_{2;4,2;4,0}(DE_2, E_4)$.

\bigskip

\subsection{Proof of Theorem 2.9:} Using any two the $f_i$'s,   
$1\le i\le 3$, in \eqref{fi} and eliminating $\Delta(z)$,  we get the following  
expression among the divisor functions. 

\begin{equation}\label{id1}
24\sum_{m=1}^{n-1}(4m^3 - 3 m^2n)\sigma(m)\sigma(n-m) = n^3 \sigma(n) - 
n^3\sigma_3(n).
\end{equation}

From Niebur's identity, we have
\begin{equation}\label{delta-e2}
\Delta(z) = E_2 D^4E_2(z) - 16 DE_2~D^3E_2 + 18 (D^2E_2)^2.
\end{equation}
It is well known that $E_2$ satisfies the following transformation 
property with respect to $SL_2({\mathbb Z})$. For $\gamma = 
\begin{pmatrix} * & *\\ c&d\\\end{pmatrix} \in SL_2({\mathbb Z})$, 
\begin{equation}
E_2(\gamma(z)) = (cz+d)^2 E_2(z) + \frac{12c(cz+d)}{2\pi i}.
\end{equation}
Using this transformation property and its successive derivatives, we see 
that the right-hand side expression in \eqref{delta-e2} is invariant under the 
stroke operation with respect to $SL_2({\mathbb Z})$ with weight $12$, provided 
the following identities are true.
\begin{eqnarray}
3(DE_2)^2 - 2E_2 D^2E_2 +2 D^3E_2 &=&0, \label{id2}\\
D^4E_2 - E_2 D^3E_2 +2 DE_2 D^2E_2 &=&0. \label{id3}
\end{eqnarray}
It is easy to see that the derivative of \eqref{id2} gives \eqref{id3} and 
so it is enough to show the identity \eqref{id2}. Note that the forms 
$(DE_2)^2$ and  $E_2 D^2E_2$  are quasimodular forms of weight $8$ and depth 
$4$. So, by using Lemma 3.3, we have 
\begin{equation*}
\begin{split}
(DE_2)^2 & = \frac{1}{5} DE_6 + 2 D^3E_2, \\
E_2 D^2E_2 & = \frac{3}{10} DE_6 + 4 D^3E_2.
\end{split}
\end{equation*}
Eliminating $DE_6$ from the above two equations, we get \eqref{id2}. 
Thus we have demonstrated another proof of the Niebur's identity. Now \eqref{id2} 
and \eqref{id3} give rise to the following identities:
\begin{eqnarray}
12\sum_{m=1}^{n-1}(5m^2 - 3 mn)\sigma(m)\sigma(n-m) &=& n^2 \sigma(n) -
n^3\sigma(n), \label{id4}\\
24\sum_{m=1}^{n-1}(3m^3 - 2 m^2n)\sigma(m)\sigma(n-m) &=& n^3 \sigma(n) -
n^4\sigma(n). \label{id5}
\end{eqnarray}
By simple manipulations of the equations \eqref{id1}, \eqref{id4} and 
\eqref{id5}, we get the identities (i) -- (iii). Identity (iv) is obtained by 
taking the difference of the identities (iii) and (iv) of Theorem 2.7. 
Next we prove (v). Eliminating the expression $\sum m^3\sigma(m)\sigma_3(n-m)$ 
from (iv) and (ix) of Theorem 2.7 and substituting for the expression  
$\sum m^2\sigma(m)\sigma_3(n-m)$ from (iv), we get (v). Eliminating $\tau(n)$ 
from (iii) and (iv) of Theorem 2.7 gives the identity (vi). To prove (vii) we 
first eliminate $\tau(n)$ between (vi) and (viii) of Theorem 2.7 and then use 
the identity (vi) of Theorem 2.9. Finally eliminating $\tau(n)$ between (ii) and 
(v) of Theorem 2.7 gives the identity (viii). This completes the proof.

\end{document}